\documentclass[11pt]{amsart}
\usepackage[english]{babel}
\usepackage{amssymb}
\usepackage{amsmath}
\newcommand{\ignore}[1]{}

\newtheorem{lemma}{Lemma}
\newtheorem{theorem}{Theorem}

\newtheorem{corollary}{Corollary}

\theoremstyle{definition}

\newtheorem{example}{Example}

\newtheorem{remark}{Remark}



\title{Norm principles for forms of higher degree permitting composition}

\author{R. W. Fitzgerald}
\email{rfitzg@math.siu.edu}
\address{
 Department of Mathematics\\
 Southern Illinois University\\
 Carbondale IL 62901-4408\\
 USA
}

\author{S. Pumpl\"un}
\email{susanne.pumpluen@nottingham.ac.uk}
\address{
 School of Mathematics\\
University of Nottingham\\
University Park\\
Nottingham NG7 2RD\\
United Kingdom
}

\begin{document}

\maketitle

\begin{abstract} Let $F$ be a field of characteristic 0 or greater than $d$.
  Scharlau's norm principle holds for  separable field extensions
 $K$ over $F$, for certain forms $\varphi$ of degree $d$ over $F$ which permit composition.
\end{abstract}

\section*{Introduction}
 Let $d\geq 2$ be an integer and let $F$ be a field of characteristic 0 or $>d$.
 Let $\varphi:V\to F$ be a  form of degree $d$ on an $F$-vector space
$V$ of dimension $n$ (i.e., after suitable identification, $\varphi$ is a homogeneous polynomial of degree $d$
in $n$ indeterminates). Let $K/F$ be a finite field extension of degree $m$.
Scharlau's norm principle (SNP) says that if $a$ is a similarity factor of $\varphi_K$, then $N_{K/F}(a)$ is a similarity
factor of $\varphi$.
Knebusch's norm principle (KNP) states that if $a$ is represented by $\varphi_F$, then
$N_{K/F}(a)$ is a product of $m$ elements represented by $\varphi$, hence lies in the subgroup of $F^\times$
generated by $D_F(\varphi)$. Both norm principles were proved  for nondegenerate quadratic forms
over fields of characteristic not 2
(cf. [Sch, II.8.6] or [L, p.~205, p.~206]). For finite extensions of semi-local regular
rings containing a field of characteristic 0,
 Knebusch's norm principle (for quadratic forms) was proved in  [Z] and for finite  \'etale extensions of
  semi-local Noetherian
domains with infinite residue fields of characteristic different from 2 in [O-P-Z].
Barquero and Merkurjev [B1,2], [B-M] generalized the norm principle to algebraic groups.

We prove Scharlau's norm principle for certain nondegenerate forms $\varphi$ of degree $d\geq 3$ which permit composition.
  Scharlau's and Knebusch's norm principle ``coincide'' for
these forms, since they permit composition in the sense of Schafer  [S] and thus satisfy
$D_K(\varphi)=G_K(\varphi)$ for all field extensions $K/F$. We explicitly compute
 the norms of some similarity factors, if $\varphi$ is the norm of a separable field extension
 or of a central simple algebra.

\section{Preliminaries}

A {\it form of degree $d$} over $F$ is a map $\varphi:V\to F$  on a finite-dimensional vector space $V$ over $F$
 such that $\varphi(a v)=a^d\varphi(v)$ for all $a\in F$, $v\in V$ and such that the map $\theta : V \times
\dots \times V \to F$ ($d$-copies) defined by
 $$\theta(v_1,\dots,v_d)=\frac{1}{d!} \sum_{1\leq i_1< \dots<i_l\leq d}(-1)^{d-l}\varphi(v_{i_1}+ \dots +v_{i_l})$$
(with $1\leq l\leq d$) is $F$-multilinear and  invariant under all permutations of its variables. The {\it dimension}
 of $\varphi$ is defined as ${\rm dim}\,\varphi={\rm dim}\, V$.
  $\varphi$ is called {\it nondegenerate}, if
$v = 0$ is the only vector such that $\theta (v, v_{2}, \dots, v_d) = 0$ for all
 $v_i \in V$. We will only study nondegenerate forms.
 Forms of degree $d$ on $V$ are in obvious one-one correspondence with homogeneous polynomials
of degree $d$ in $n={\rm dim}\, V$ variables.
 If $\varphi$ is represented by $a_1x_1^d+\ldots +a_mx_m^d$ ($a_i\in F^\times$), we use the notation $\varphi=\langle  a_1,\ldots
,a_n\rangle  $ and call $\varphi$ {\it diagonal}.

Two forms $(V_i,\varphi_i)$ of degree $d$, $i=1,2$, are called {\it isomorphic} (written
$(V_1,\varphi_1)\cong (V_2,\varphi_2)$ or just $\varphi_1\cong\varphi_2$) if there exists a bijective linear map
$f:V_1\to V_2$ such that $\varphi_2(f(v))=\varphi_1(v)$ for all $v\in V_1.$

 Let $(V,\varphi)$ be a form over $F$ of degree $d$ in $n$ variables over $F$.
 An element $a\in F$ is {\it represented by} $\varphi$ if there is an $v\in V$ such that $\varphi(v)=a$.
 An element $a\in F^{\times}$ such that $\varphi\cong a\varphi$ is called a {\it similarity factor} of $\varphi$.
Write $D_F(\varphi)= \{a \in F^\times\mid\varphi(x)=a \text{ for some } x\in V\}$ for the set of non-zero
elements represented by $\varphi$ over $F$ and $G_F(\varphi)=\{a\in F^\times\mid\varphi \cong a\varphi\}$ for the group of similarity factors of $\varphi$ over $F$.
The subscript $F$ is omitted if it is
clear from the context that $\varphi$ is a form over the base field $F$.
$\varphi$  is called \emph{round} if $D(\varphi)\subset G(\varphi)$.

 A nondegenerate form $\varphi(x_1,\ldots,x_n)$ of degree $d$ in $n$ variables
\emph{permits composition} if $\varphi(x)\varphi(y)=\varphi(z)$ where $x$, $y$ are systems of $n$
indeterminates and where each $z_l$ is a bilinear form in $x,y$ with coefficients in $F$.
 In this case the vector space
$V=F^n$ admits a bilinear map $V\times V\rightarrow V$ which can be viewed as the multiplicative structure of a
nonassociative $F$-algebra and $\varphi(vw)=\varphi(v)\varphi(w)$ holds for all $v, w\in V$.
Note that the form $\varphi$ here is nondegenerate if and only if the underlying (automatically alternative)
$F$-algebra is separable (Schafer [S]).
For instance, every norm of a central simple algebra or of a separable finite field extension over $F$ is nondegenerate
and permits composition.

\begin{remark} (i) There are two types of forms $\varphi$ of degree $d$ over $F$ for which SNP trivially holds:\\
(a) if $G_F(\varphi)=F^{\times }$;\\
(b) if $G_K(\varphi)=K^{\times d}$ for every field extension $K$ over $F$.\\
 (ii) Let $\varphi$ be a diagonal form over $F$ of degree $d\geq 3$. If ${\rm dim}\,\varphi=1$ or ${\rm dim}\,\varphi\in
\{sd+1,sd-1\}$ for some integer $s\geq 1$, then $G_K(\varphi)=K^{\times d}$ for every finite field extension $K$ over
 $F$ [Pu, Proposition 1 (i)].
 Hence $\varphi$ trivially satisfies SNP for all field extensions $K$ over $F$ by (i).
 Moreover, every form $\langle a,a,\dots,a\rangle$ of degree $d\geq 3$ satisfies
 $G_K(\varphi)=K^{\times d}$ for all field extensions $K$ over $F$ [Pu, Lemma 9 (ii)], hence SNP.
\\
(iii) If $\varphi$ is the determinant of the $d$-by-$d$ matrices over $F$, then $G_K(\varphi)=K^{\times }$ for all field extensions $K$ over $F$,
 hence SNP holds for all field extensions of $F$ by (i).\\
(iv) The cubic norm $\varphi$ of a reduced Freudenthal algebra  $J=H_3(C,\Gamma)$ [KMRT,  p.~516],
$C$ a composition algebra over $F$ or 0, trivially satisfies SNP for all field extensions $K$ of $F$,
 because $D_K(\varphi)=G_K(\varphi)=K^\times$.\\
(v) Suppose the base field $F$ has characteristic 0 or greater than $d+1$.  Let $\varphi_0:V\to F$ be a form of degree $d$, then the form $\varphi(a+u)=a\varphi_0(u)$, $a\in F$, $u\in V$
of degree $d+1$ satisfies $G_K(\varphi)=D_K(\varphi)=K^\times$ for all field extensions $K$ over $F$, hence SNP.
\end{remark}

\section{Forms satisfying Scharlau's norm principle}

\subsection{Norms of field extensions}

\begin{theorem}
Let $F$ be a field of arbitrary characteristic (that is, we drop our standing assumptions
 on $\text{char}(F)$). Let $\widetilde{F}=F(\alpha)$ be a  field
extension of $F$ of degree $d$ and $\varphi=N_{\widetilde{F}/F}$.
 Suppose that $K/F$ is a finite field extension which is linearly disjoint to $\widetilde{F}$ over $F$.
If $e\in K^\times$ is represented by $\varphi_K$, then $N_{K/F}(e)$ is represented by $\varphi$ and thus
 $$N_{K/F}(G_K(\varphi_K))\subset G_F(\varphi).$$
\end{theorem}

\begin{proof} For all field extensions $L/F$, $D_L(\varphi_L)=G_L(\varphi_L)$. We have
$K(\alpha)\cong F(\alpha)\otimes_FK $ and $\varphi_K = N_{K(\alpha)/K} $.
If $e\varphi_K\cong\varphi_K$, then $e=\varphi_K(z_1+\alpha z_2+\dots+\alpha^{d-1}z_d)$
  with $z_i\in K$ and
$$\begin{array}{l}
N_{K/F}(\varphi_K (z_1+\alpha z_2+\dots+\alpha^{d-1}z_d))=\\
N_{K/F}(N_{K(\alpha)/K}(z_1+\alpha z_2+\dots+\alpha^{d-1}z_d))=\\
N_{F(\alpha)/F}(N_{K(\alpha)/F(\alpha)}(z_1+\alpha z_2+\dots+\alpha^{d-1}z_d))=\\
N_{F(\alpha)/F}( a_1+\alpha a_2+\dots+\alpha^{d-1}a_d)=\\
\varphi ( a_1+\alpha a_2+\dots+\alpha^{d-1}a_d)\in G_F(\varphi)
\end{array}$$
for suitable $a_i\in F$.
\end{proof}

 This simple trick which even gives an explicit identity for $N_{K/F}(e)$ in terms of the $a_i$'s,
 was used in [F] to compute norms for the quadratic form $\langle 1,1\rangle$.

We cannot omit the condition that $K$ and  $F(\alpha)$ are linearly disjoint in Theorem 1, since it is needed in
the proof to guarantee that $\varphi_K = N_{K(\alpha)/K} $ (see [La, p.~267] for an example of two extensions
which are not linearly disjoint, where this does not work). However, in
Theorem 4 we will show that for $\varphi=N_{\widetilde{F}/F}$, SNP indeed holds  for all
 Galois extensions $K/F$.

Let $\bar F$ be the algebraic closure of $F$ and $K/F$ be a finite separable
 field extension of degree $n$ in $\bar F$. Put  $\widetilde{K}=K(\alpha)$.
The restrictions of the  $\widetilde{F}$-monomorphisms  $\sigma_1,\dots,\sigma_n$ of $\widetilde{K}$ into $\bar F$
 give the $F$-monomorphisms of $F$ into $\bar F$.

\begin{theorem} Let $\widetilde{F}=F(\alpha)$ be a
field extension of $F$ of degree $d$ such that $\alpha^d=c$ and let $\varphi=N_{\widetilde{F}/F}$ be its norm.
Let $K/F$ be a quadratic field extension and
let $\sigma :K \to K$ be the non-identity $F$-automorphism. Let $z_1,\dots,z_n\in K$. Then
$$\begin{array}{l}
N_{K/F}(\varphi(z_1,z_2,\dots,z_n))=\\
\varphi(A_0+cA_1,A_2+cA_3,A_4+cA_5,\dots, A_{2d-4}+cA_{2d-3}, A_{2d-2})
\in G_F(\varphi)
\end{array}$$
where $A_i\in F$ for $0\leq i\leq 2d-2$.
\end{theorem}

\begin{proof} First suppose $[K(\alpha ) :K]=d$. Set $z=\sum_{i=1}^d z_i\alpha^{i-1}\in K(\alpha )$.
Then we have that $\varphi_K(z_1, z_2, \ldots , z_d)=N_{K(\alpha )/K}(z)$ and so
$$
N_{K/F}(\varphi_K(z_1, z_2, \ldots ,z_d))=N_{K(\alpha )/K}(z)=N_{F(\alpha )/F}(N_{K(\alpha )/K}(z)).
$$
 Let $id$, $\sigma$ be the $\widetilde{F}$-monomorphisms of
$\widetilde{K}=K(\alpha)$ into $\bar F$. If $0\not=z$, then
$$\begin{array}{l}
N_{K(\alpha)/F(\alpha)}(z)=(z_1+\alpha z_2+\dots+\alpha^{d-1}z_d) \sigma (z_1+\alpha z_2+\dots+\alpha^{d-1}z_d)=\\
z_1\sigma (z_1)+c( z_2\sigma ( z_d)+ z_3\sigma ( z_{d-1})+\dots+ z_d\sigma ( z_{2}))+\\
\alpha [z_1\sigma (z_2)+z_2\sigma (z_1)+c( z_3\sigma ( z_d)+ z_4\sigma ( z_{d-1})+\dots+ z_d\sigma ( z_{3}))]+\\
\alpha^2 [z_1\sigma (z_3)+z_2\sigma (z_2)+z_3\sigma (z_1)+c( z_4\sigma ( z_{d})+\dots+ z_d\sigma ( z_{4}))]+\\
\vdots\\
\alpha^{d-1} [z_1\sigma (z_d)+z_2\sigma (z_{d-1})+\dots+ z_d\sigma ( z_{1})]=\\
A_0+cA_1+\alpha [A_2+cA_3]+\alpha^2 [A_4+cA_5]+\dots+\alpha^{d-2} [A_{2d-4}+cA_{2d-3}]+\alpha^{d-1}A_{2d-2}.
\end{array}$$
with the $A_i$'s given by
$$\begin{array}{l}
A_0=z_1\sigma (z_1)\\
A_1= z_2\sigma ( z_d)+ z_3\sigma ( z_{d-1})+\dots+  z_{d-1}\sigma ( z_{3})+z_d\sigma ( z_{2})\\
A_2=z_1\sigma (z_2)+z_2\sigma (z_1)\\
A_3=z_3\sigma ( z_d)+ z_4\sigma ( z_{d-1})+\dots+ z_{d-1}\sigma ( z_{4})+z_d\sigma ( z_{3})\\
A_4=z_1\sigma (z_3)+z_2\sigma (z_2)+z_3\sigma (z_1)\\
A_5= z_4\sigma ( z_{d})+\dots+z_d\sigma ( z_{4})\\
A_6= z_1\sigma ( z_{4})+ z_2\sigma ( z_{3})+ z_3\sigma ( z_{2})+z_4\sigma ( z_{1})\\
A_7= z_5\sigma ( z_{d})+\dots+z_d\sigma ( z_{5})\\
\vdots\\
A_{2d-4}=z_1\sigma (z_{d-1})+z_2\sigma (z_{d-2})+z_3\sigma (z_{d-3})+\dots+ z_{d-1}\sigma ( z_{1})\\
A_{2d-3}=z_d\sigma (z_d)\\
A_{2d-2}=z_1\sigma (z_d)+z_2\sigma (z_{d-1})+\dots+ z_{d-1}\sigma ( z_{2})+z_d\sigma ( z_{1}).
\end{array}$$
 $A_i\in F$  for all $i$, since each $A_i$ is invariant under $\sigma$.
 Applying $N_{F(\alpha )/F}$ yields the assertion.
 The result is a polynomial identity over $K$ and so it also holds when $[K(\alpha ) : K] < d$.
\end{proof}

A similar argument should work for arbitrary  field extensions $\widetilde{F}/F$ of degree $d$, so that one should be
able to prove  Theorem 2 more generally. However, the computations
become extremely tedious. We look at  one more special case which  is proved analogously:

\begin{theorem} Let $\widetilde{F}=F(\alpha)$ be a
field extension of $F$ of degree $d$ such that $\alpha^d=b\alpha+c$ and $\varphi=N_{\widetilde{F}/F}$ its norm.
Let $K/F$ be a quadratic field extension and let $\sigma :K \to K$ be the non-identity $F$-automorphism.
 Let $z_1,\dots,z_n\in K$.
$$\begin{array}{l}
N_{K/F}(\varphi_K(z_1,\dots,z_n))=\\
\varphi(A_0+cA_1,A_2+cA_3+bA_1,A_4+cA_5+bA_3,A_6+cA_7+bA_5,\dots,\\
 A_{2d-4}+cA_{2d-3}+bA_{2d-5},A_{2d-2}+bA_{2d-3})\in G_F(\varphi),
\end{array}$$
with $A_i\in F$ for $0\leq i\leq 2d-3$  defined as in Theorem 2.
\end{theorem}

\begin{corollary} If $\varphi$ is the norm of a field extension of degree $d$ of $F$ of the  kind given in
Theorem 2 or 3, then $\varphi$ satisfies SNP for all quadratic field extensions of $F$.
\end{corollary}

We state the following

\smallskip\noindent
{\bf Conjecture.} {\it Let $\varphi$ be the norm of a  field extension of degree $d$ of $F$. Then
 SNP holds for $\varphi$ for all separable field extensions $K/F$.}

\smallskip
It is highly likely the above method of proof can be used to prove this.
However, it is not clear how to explicitly do the above calculations in this general case. If we
restrict ourselves to forms of prime degree and Galois extensions, however, we can prove SNP because of the following easy observation:

\begin{remark} (i) Let $\varphi$ be a form of degree $d$ over $F$. Let $K/F$ be a finite field extension.
Suppose we have $a\varphi_K\cong\varphi_K$ for some $a\in K^\times$.\\
(a) If $[K:F(a)]=dm$ then a straightforward calculation shows that $N_{K/F}(a)\in F^{\times d}\subset G(\varphi)$.\\
(b) If $a\in F$ then trivially $N_{K/F}(a)\in F^{\times d}\subset G(\varphi)$.\\
(ii) Let $\varphi$ be a form of prime degree $p$ over $F$. Then SNP holds for $\varphi$ for all field extensions
 of degree $p^r$ for some integer $r>0$ by (a).
\end{remark}

\begin{corollary}  Let $\widetilde{F}=F(\alpha)$ be a  field extension of $F$ of prime degree $p$ and
$\varphi=N_{\widetilde{F}/F}$ its norm. Then SNP holds for all extensions $K/F$ which can be written as a tower of fields
$$F=K_0\subset K_1\subset\dots\subset K_n=K$$
such that $[K_{i+1}:K_i]$ either is a power of $p$  or coprime with $p$.
\end{corollary}

\begin{proof} Let $a\varphi_K\cong\varphi_K$
for $a\in K^\times$.
 If $[K:K_{n-1}]=p^r$ then $N_{K/K_{n-1}}(a)\varphi_{K_{n-1}}\cong\varphi_{K_{n-1}}$ by Remark 2 (ii).
 If $[K:K_{n-1}]$ is coprime to $p$ then $N_{K/K_{n-1}}(a)\varphi_{K_{n-1}}\cong\varphi_{K_{n-1}}$ by Theorem 1.
Since
$N_{K/F}= N_{K_{1}/F} \circ  \dots \circ N_{K_{n-1}/K_{n-2}}\circ N_{K/K_{n-1}}$, a straightforward induction
 implies the assertion.
\end{proof}

\begin{theorem} Let $\widetilde{F}=F(\alpha)$ be a  field
extension of $F$ of prime degree $p$ and $\varphi=N_{\widetilde{F}/F}$ its norm. Then SNP
holds  for all finite Galois  extensions $K/F$.
\end{theorem}

\begin{proof}  Suppose that $K/F$ is a Galois extension and $[K:F]=p^rm$ where $m$ is not divisible by $p$.
There is a separable intermediate field $K_0$ of $K/F$ such that $[K:K_0]=p^r$ and $[K_0:F]=m$.
Apply Corollary 2.
\end{proof}

\subsection{Norms of central simple algebras}

We now turn to the (reduced) norm forms of central simple algebras over $F$.
 Let $\varphi=N_{A/F}$ be the norm of a central simple algebra $A$ of degree $d$ over $F$. Then
  SNP holds for all finite separable field extension [B-M, 3.1].
For the split central simple algebra $A\cong{\rm Mat}_d(F)$,
 $\varphi$ trivially satisfies SNP for all field extensions of $F$.

 If $A$ is a division algebra then SNP holds for all finite field extensions:

 Let $K/F$ be a finite field extension of degree $n$.
For $\alpha\in F$, $\rho_{\alpha} : K\to K$, $\rho_{\alpha}(x)=\alpha x$ is left multiplication with $\alpha$.
Fix a basis $B =\{ w_1, w_2, \ldots ,w_n\}$ of $K/F$. Let $\rho (\alpha )$ be the matrix representation of
$\rho_{\alpha}$ with respect to $B$. The map $\rho : K\to M_n(F)$ is an injective ring homomorphism and
the norm is given by $N_{K/F}(\alpha)=\det \rho (\alpha )$.

Let $A$ be a central simple algebra over $F$.
Pick $\Delta =\sum_{i=1}^n \alpha_iw_i$, where $\alpha_i\in A$ and so $\Delta\in \bar A =A\otimes{K}$.
Again, $\rho_{\Delta} : \bar A\to \bar A$ is left multiplication and $\rho (\Delta)$ is the matrix, with entries in $A$, of
$\rho_{\Delta}$ with respect to $B$.
For the proof of the next theorem we need the following observation:

\begin{lemma}
$\rho (\Delta )=\sum_{i=1}^n \alpha_i\rho (w_i)$.
\end{lemma}

\begin{proof} Let $a\in \bar A$. Then $\rho_{\Delta}(a)=\sum \alpha_iw_ia =\sum \alpha_i \rho_{w_i}(a)$. Hence $\rho_{\Delta}=\sum \alpha_i\rho_{w_i}$ and for matrices $\rho (\Delta )=\sum \alpha_i \rho (w_i)$.
\end{proof}

Let $A$ be a central simple division algebra over $F$ with basis $\epsilon_1,\dots,\epsilon_m$.
Let $A^\times$ be the invertible elements in $A$ and $C(D^\times)=[A^\times,A^\times]$ be the commutator
subgroup. Put $\bar A=A\otimes_F K$. Let ${\rm det}:{\rm GL}_n A\to A^\times/C(A^\times)$ be the Dieudonn\'{e}
determinant.
There is a polynomial
$G\in F[x_1,\dots,x_m]$ such that for any extension $L/F$ the norm from $A\otimes L\to L$ is given by
$$N(\sum_{i=1}^m l_i\epsilon_i)=G(l_1,\dots,l_m).$$
We write
$$G(*\, l_k *) \text{ for } G(l_1,\dots,l_k,\dots, l_m).$$

\begin{theorem}
Let $A$ be a central simple division algebra over $F$. Let $K/F$ be a finite extension (which need not be separable). Then
$$
N_{K/F}(N_{\bar A/K}(\Delta )=N_{A/F}(\det \rho (\Delta )).
$$
\end{theorem}

\begin{proof} The matrices $\rho (w_1), \rho (w_2), \ldots , \rho (w_n)$ commute and so have a common eigenvector. A simple induction argument shows that there is a matrix $P$, over the algebraic closure $\bar F$, such that each $P^{-1}\rho (w_i)P$ is upper triangular. Let the diagonal entries of $P^{-1}\rho (w_i)P$ be denoted by $d_{ij}$, $1\le j\le n$.

We compute both sides starting with the right-hand side: By  Lemma 1,
$$
P^{-1}\rho (\Delta )P =\begin{pmatrix} \sum_i \alpha_id_{i1} & & &\\ & \sum_i \alpha_id_{i2} & & *\\0 & & \ddots & \\ & & & \sum_i \alpha_id_{in}\end{pmatrix}.
$$
Now Dieudonn\'e's determinant [P, p.~308] satisfies $\det (P^{-1}AP)=\det A$ and the determinant of an upper triangular matrix is the product of the diagonal elements
 (in [A, p. 163], the first is consequence h), the second follows from [A, Theorem 4.4]). Hence
$$
\det\rho (\Delta )=\prod_{j=1}^n \left(\sum_{i=1}^n \alpha_id_{ij}\right).
$$

Write $\alpha_i =\sum_{k=1}^m a_{ik}\epsilon_k$ where $ a_{ik}\in F$.
For the right-hand side we know that
\begin{eqnarray*}
\det \rho (\Delta ) &=& \prod_{j=1}^n \sum_{k=1}^m \left(\sum_{i=1}^n a_{ik}d_{ij}\right)\epsilon_k,\\
N_{A/F}(\det \rho (\Delta )) &=& \prod_{j=1}^n G(* \sum_{i=1}^n a_{ik}d_{ij} \,*).
\end{eqnarray*}
For the left-hand side we have
\begin{eqnarray*}
\Delta &=& \sum_{k=1}^m\left(\sum_{i=1}^n a_{ik}w_i\right)\epsilon_k,\\
N_{\bar A/K}(\Delta ) &=& G(* \sum_{i=1}^n a_{ik}w_i\, *).
\end{eqnarray*}
As $\rho$ is a ring homomorphism, $\rho(G(*\, u_k*))=G(*\, \rho(u_k)*)$. Thus
$$
N_{K/F}(N_{\bar A/K}(\Delta )) = \det G(*\, \sum_{i=1}^n a_{ik}\rho (w_i)\,*).
$$
Conjugation by $P$ is also a ring homomorphism, so
$$
N_{K/F}(N_{\bar A/K}(\Delta )) = \det G(* \, \sum_{i=1}^n a_{ik}P^{-1}\rho (w_i)P\, *).$$
We conclude that $G(*\, \sum_{i=1}^n a_{ik}P^{-1}\rho (\beta )^iP \,*)=$
\begin{eqnarray*}
 && G\left(*\quad  \begin{pmatrix} \sum_i a_{ik}d_{i1} & & &\\ & \sum_i a_{ik}d_{i2} & & *\\0 & & \ddots & \\ & & & \sum_i a_{ik}d_{in}\end{pmatrix}\quad
*\right)=\\
 && \begin{pmatrix} G(*\, \sum_i a_{ik}d_{i1}\, *) & & &\\ & G(*\, \sum_i a_{ik}d_{i2}\, *) & & *\\0 & & \ddots & \\ & & & G(*\, \sum_i a_{ik}d_{in}\, *)\end{pmatrix}.
\end{eqnarray*}
Hence
$$
N_{K/F}(N_{\bar A/K}(\Delta )) = \prod_{j=1}^n G(*\, \sum_{i=1}^n a_{ik}d_{ij}\,*),
$$
the same as the right-hand side, proving the identity.
\end{proof}

\begin{theorem} Let $\varphi$ be the norm of a central simple division algebra $A$ over $F$. Then SNP holds
for all finite field extensions of $F$.
\end{theorem}

\begin{proof} The proof is analogous to the one given in [F, Lemma 2.1] for the norms of a quaternion
division algebra: Let $\epsilon_1,\dots,\epsilon_m$ be a basis for $A$ as a $F$-vector space
(where $m=d^2$ if $d$ is the degree of $A$).
For $z_i\in K$ and $z=\epsilon_1z_1+\epsilon_2z_2+\dots+\epsilon_{m}z_{m}$, we have
$$\begin{array}{l}
n_{K/F}(\varphi_K(z))=\\
n_{K/F}(n_{\widetilde{D}/K}(z))=\\
n_{A/F}({\rm det}(\rho(z)))=\\
n_{A/F}(\epsilon_1a_1+\epsilon_2a_2+\dots+\epsilon_{m}a_{m})
\end{array}$$
for suitable $a_i\in F$. (The second equality holds by Theorem 5.)
\end{proof}

\begin{corollary}  Let $\varphi$ be the norm of a central simple algebra $A$ over $F$ of prime degree. Then SNP holds
for all finite field extensions of $F$.
\end{corollary}

\begin{remark} Let $K=F(\sqrt{c})$ be a quadratic field extension and $A$ a division algebra over $F$ of degree $d$.
Let $z_i=u_i+v_i\sqrt{c}\in K$ and $z=z_1\epsilon_1+z_2\epsilon_2+\dots+z_{d^2}\epsilon_{d^2}$, then $z=x+y\sqrt{c}$ with
$x=u_1\epsilon_1+u_2\epsilon_2+\dots+u_{d^2}\epsilon_{d^2}$ and $y=v_1\epsilon_1+v_2\epsilon_2+\dots+v_{d^2}\epsilon_{d^2}$.
We obtain, more explicitly than above (similar as in [F, 2.2]):
$$\begin{array}{l}
N_{K/F}(\varphi_K(z))=
N_{A/F}({\rm det}(\rho(z)))=
N_{A/F}(y(xy^{-1}x-dy))\in D_F(N_{A/F}).
\end{array}$$
In particular, if $A$ has degree 3, then we can also write
$$\begin{array}{l}
N_{K/F}(\varphi_K(z))=
\frac{1}{N_{A/F}(y)}N_{A/F}(xy^\sharp x-dN_{A/F}(y)y)
\end{array}$$
with $x^\sharp=x^2-T_{A/F}(x)x+S_{A/F}(x)1_A$ [KMRT, p.~470].
\end{remark}

\begin{remark} Let $F$ be an infinite field. Suppose there are $f,g\in F[X_1,\dots ,X_n]$ such that
$$f(u_1,\dots,u_n)^m=g(u_1,\dots,u_n)^m \text{ for all }u_i\in F.$$
Then there is an $m$th root of unity $\mu$ in $F$ such that
$$f(u_1,\dots,u_n)=\mu g(u_1,\dots,u_n)  \text{ for all }u_i\in F.$$
\end{remark}

\begin{lemma} Let $\varphi_1:V\to F$ be a form over $F$ of degree $d_1$ which satisfies SNP.
 Put $\varphi(u)=\varphi_1(u)^m$ for some integer $m\geq 2$.\\
(i) If $D_K(\varphi_1)=G_K(\varphi_1)$ for all finite field extensions $K/F$, then $\varphi$ satisfies SNP
 for all finite field extensions.\\
(ii) If $F$ is an infinite field, then  $\varphi$ satisfies SNP  for all finite field extensions.
\end{lemma}

\begin{proof} (i) Let $a\varphi_K\cong\varphi_K$. Then $a=b^m$ for a suitable $b=\varphi_1(w)$ and therefore
$N_{K/F}(b)\in G_F(\varphi_1)$ for every finite field extension $K/F$. This implies $N_{K/F}(a)=N_{K/F}(b^m)\in G_F(\varphi)
=D_F(\varphi)$.\\
(ii)  Let $a\varphi_K\cong\varphi_K$ for some
finite field extension $K/F$. Then there is an isomorphism $f:V\to V$ such that
$a\varphi_{1,K}(x)^m=\varphi_{1,K}(f(x))^m$ for all $x\in V$. Let $x$ be an anisotropic vector in $V$, then
$a=(\varphi_{1,K}(f(x))/\varphi_{1,K}(x))^m$ is an $m$th power in $F$, hence write $a=b^m$ for some $b\in K^\times$.
From $b^m\varphi_{1,K}^m\cong\varphi_{1,K}^m$ we conclude that $\mu b\varphi_{1,K}\cong\varphi_{1,K}$
for some $m$th root of unity $\mu$ in $F$ (Remark 5). Therefore $N_{K/F}(\mu b)\in G_F(\varphi_1)$ for every finite field extension $K/F$ and thus
also $N_{K/F}(\mu b)^m=N_{K/F}(a)\in G_F(\varphi)$.
\end{proof}

\begin{lemma} (i) Let $\varphi_i:V_i\to F$ be two forms over $F$ of degree $d_i$ which satisfy SNP
for all finite field extensions $K/F$.
 Put $\varphi:V_1\oplus V_2\to k,$ $\varphi(u)=\varphi_1(u_1)\varphi_2(u_2)$ for $u=u_1+u_2$, $u_i\in V_i$.
  If $D_K(\varphi_i)=G_K(\varphi_i)$ for all finite field extensions $K/F$, then $\varphi$ satisfies SNP
   for all finite field extensions.\\
  (ii) Let $F'/F$ be a finite separable field extension and $\varphi_0:V\to F'$ be a form over $F'$.
Let $\varphi=N_{F'/F}(\varphi_0)$. Suppose that $(\varphi_0)_{L'}$ is a round form for all
finite field extensions $L'$ of $F'$ and that SNP holds for $\varphi_0$ for all
finite field extensions $L'$ of $F'$. Then $\varphi=N_{F'/F}(\varphi_0)$ satisfies SNP for all finite field extensions
 $K$ of $F$ which are linearly disjoint with $F'$ over $F$.
\end{lemma}

\begin{proof} (i) By [Pu], $\varphi_K$ is a round form. Let $a\varphi_K\cong\varphi_K$.
Then $a=\varphi_{1,K}(w_1)\varphi_{2,K}(w_2)$ and by assumption, $N_{K/F}(\varphi_{i,K}(w_i))\in G_F(\varphi_i)$
for $i=1,2$. This immediately yields
$N_{K/F}(\varphi_{1}(w_1))N_{K/F}(\varphi_{2}(w_2))=N_{K/F}(\varphi_{1}(w_1)\varphi_{2}(w_2))=
N_{K/F}(a)\in G_F(\varphi)$.\\
 (ii)  Let $K$ be a finite field extension of $F$ which is linearly disjoint with $F'$ over $F$.
 Then
 $$\varphi_{K}=N_{K'/K}((\varphi_0)_{K'})$$
with $K'=F'\cdot K$ the composite of $F'$ and $K$ (i.e., the homogeneous polynomials defining the forms are equal).
Since $(\varphi_0)_{K'}$ is round by assumption, $D_{K'}((\varphi_{0})_{K'} )=G_{K'}((\varphi_{0})_{K'})$, and
 $\varphi_K$ is a round form by [Pu].\\
Let $a\varphi_K\cong\varphi_K$. Since $\varphi_K$ is round, $a=N_{K'/K}((\varphi_{0})_{K'}(z_0))$ for some $z_0\in K'$.
As $(\varphi_0)_{K'}$ is round, we have
$$(*)\,\,\,\,\,\,\,\,((\varphi_0)_{K'}(z_0))(\varphi_0)_{K'}\cong (\varphi_0)_{K'}.$$
$\varphi_{0}$ satisfies SNP for all field extensions of $F'$ by assumption, hence
$$N_{K'/F'}((\varphi_0)_{K'}(z_0))\varphi_{0}\cong\varphi_{0}$$
 and so
$N_{F'/F}(N_{K'/F'}((\varphi_0)_{K'}(z_0)))\varphi\cong \varphi.$
Hence
$$N_{F'/F}(N_{K'/F'}((\varphi_0)_{K'}(z_0)))=N_{K/F}(N_{K'/K}((\varphi_0)_{K'}(z_0)))=
N_{K/F}(a)\in G_F(\varphi).$$
\end{proof}

Similarly, we obtain:
\begin{theorem} (a) Let $L/F$ be a separable field extension such that there is an intermediate field
$F\subset F'\subset L$ with $[L:F']=p$ prime and $[F':F]$ a Galois extension. Let $\varphi=N_{L/F}$ be its norm.
Then SNP holds for all Galois extensions $K/F$ of degree $[K:F]$ coprime to $[F':F]$.\\
(b) Let $F'/F$ be a finite separable field extension and $\varphi_0:V\to F'$ be a form over $F'$ of prime degree $p$.
Let $\varphi=N_{F'/F}(\varphi_0)$. Suppose that $(\varphi_0)_{L'}$ is a round form for all
finite field extensions $L'$ of $F'$. Then $\varphi=N_{F'/F}(\varphi_0)$ satisfies SNP for all field extensions
 $K$ of $F$ of degree $p^r$ coprime to $[F':F]$.
\end{theorem}

\begin{proof} (a) Let $\varphi_0=N_{K/F'}: L\to F'$, which has prime degree $p$. $(\varphi_0)_{L'}$
is a round form for all finite field extensions $L'$ of $F'$.

Let $K/F$ be a Galois extension of degree coprime to $[F' : F]$. Then $K'=F'\cdot K$ is linearly disjoint from
$F'$ over $F$. The proof of Lemma 3 (ii) holds up to (*). By Theorem 4, SNP holds for $\varphi_0$ for all Galois
extensions of $F'$, in particular, for $K'$. So (*) yields
$$
N_{K//F'}((\varphi_0)_{K'}(z_0))\varphi_0\cong \varphi_0,
$$
and $N_{K/F}(a)\in G_F(\varphi )$.
\\
(b) Let $K$ be a field extension of degree $p^r$ which is coprime to $[F' : F]$
 and set $K'=F'\cdot K$. Then $[K' : F']=p^r$ and $K'$ is linearly disjoint from $F'$ over $F$.
 Again the proof of Lemma 3 (ii) holds up to (*). By Remark 2 (ii), SNP holds for $\varphi_0$ for all extensions
$K/F'$ of degree a power of $p$, in particular, for $K'$. So (*) yields $N_{K/F}(a)\in G_F(\varphi )$ as before.
\end{proof}

Theorem 1, applied to the setting in (a), only proves SNP for all extensions of $F$ of degree coprime to $p\cdot [F':F]$.

Forms $\varphi_0$ over $F'$ which satisfy the conditions of Theorem 7 (b) are not only those permitting composition, but also
forms permitting Jordan composition of prime degree over fields of characteristic 0 or greater than $2d$,
 e.g. the reduced cubic norm of an Albert algebra.
\begin{example} Let $\varphi_0=\langle\langle a_1,\dots,a_r\rangle\rangle$ ($a_i\in F^\times$) be an
anisotropic $r$-fold quadratic Pfister form.
If $K=F(\sqrt{c})$ is a quadratic field extension, then
$$\begin{array}{l}
N_{K/F}(\varphi_0)(u_1,w_1,\dots,u_{2^r},w_{2^r})=\\
(\langle\langle a_1,\dots,a_r,c\rangle\rangle)^2 (u_1,u_2,\dots,u_{2^r},w_1,w_2,\dots,w_{2^r})-4c\varphi_0(u_1w_1,\dots,
u_{2^r}w_{2^r})
\end{array}$$
is an anisotropic quartic form of dimension $2^{r+1}$ which satisfies SNP  for all finite field extensions of $F$
which are linearly disjoint with $K$ over $F$.

If $F$ contains a primitive third root of unity and $K=F(\sqrt[3]{c})$ is a cubic Kummer field extension, then
$$\begin{array}{l}
N_{K/F}(\varphi_0)(u_1,v_1,w_1,\dots,u_{2^r},v_{2^r},w_{2^r})=\\
(\langle\langle a_1,\dots,a_r,2c\rangle\rangle)^3 (u_1,\dots,u_{2^r},v_1w_1,\dots,v_{2^r}w_{2^r})\\
+c(c\langle\langle a_1,\dots,a_r\rangle\rangle\perp 2\langle\langle a_1,\dots,a_r\rangle\rangle)^3
(w_1,\dots,w_{2^r},u_1v_1,\dots,u_{2^r}v_{2^r})\\
+c^2(\langle\langle a_1,\dots,a_r\rangle\rangle\perp 2\langle\langle a_1,\dots,a_r\rangle\rangle)^3
(v_1,\dots,v_{2^r},u_1w_1,\dots,u_{2^r}w_{2^r})\\
-3c[(\langle\langle a_1,\dots,a_r,2c\rangle\rangle(u_1,u_2,\dots,u_{2^r},v_1w_1,\dots,v_{2^r}w_{2^r}))\\
\cdot
((c\langle\langle a_1,\dots,a_r\rangle\rangle\perp 2\langle\langle a_1,\dots,a_r\rangle\rangle)
(w_1,\dots,w_{2^r},u_1v_1,\dots,u_{2^r}v_{2^r}))\\
\cdot
(\langle\langle a_1,\dots,a_r\rangle\rangle\perp 2\langle\langle a_1,\dots,a_r\rangle\rangle)
(v_1,\dots,v_{2^r},u_1w_1,\dots,u_{2^r}w_{2^r}))]
\end{array}$$
is an anisotropic form of degree 6 and dimension $3\cdot 2^{r}$ which satisfies SNP
 for all finite field extensions of $F$ which are linearly disjoint with $K$ over $F$.
\end{example}

 There exists a nondegenerate form $\varphi$ of degree $d>2$ permitting composition on  a finite dimensional unital $F$-algebra
  $A$ if and only if
$A$ is a separable alternative algebra and $\varphi$ is one of the following forms, for some
integers $s_1,\ldots , s_r>0$: write $A$ as direct sum of simple ideals
$A=A_1\oplus\dots\oplus A_r$
 with the center of each $A_i$ a separable field extension $F_i$ of $F$. Any $a\in A$
can be written uniquely as $a=a_1+\ldots +a_r,\, a_i\in A_i$ and any nondegenerate form $\varphi$ on $A$ permitting
composition can be written as
$$\varphi(a)=N_1(a_1)^{s_1}\cdots N_r(a_r)^{s_r},$$
 where $d=d_1s_1+\ldots +d_rs_r$, and where $N_i$ is the generic norm of the $F$-algebra $A_i$ of degree $d_i$
 [S]. If SNP holds for all $N_i$ then it holds for $\varphi$ (Lemma 2, 3).

\begin{theorem} If $\varphi$ is a nondegenerate cubic form  over $F$ which
permits composition, which is not the norm of a cubic field extension, then SNP holds for all finite separable
field extensions of $F$.
\\ If $\varphi=N_{F'/F}$ is the norm of a cubic field extension, SNP holds  for every finite Galois extension of $F$,
and for every field extension $K/F$ of degree either a power of 3 or of degree coprime to 3.
\end{theorem}

It remains open for now whether or not SNP holds in the remaining cases.

\begin{proof} We have either $\varphi\cong\langle  1\rangle$, $\varphi$ is the norm of a cubic field extension,
 of a central simple $F$-algebra of degree $3$ or $\varphi(a+x)=a N_C(x)$ for $a\in F,
\, x\in C$, $C$ a composition algebra over $F$. SNP holds for all separable finite field
extensions of $F$, unless $\varphi$ is the norm of a cubic field extension, in which case
Theorems 1 and 4 yield
 the assertion.
\end{proof}

\begin{remark} (i) If the conjecture is true, we have SNP (for all separable field extensions $K/F$) for all nondegenerate
cubic forms permitting composition. So far, we have proved it for all Galois extensions $K/F$,  for all nondegenerate forms
of degree 3  permitting composition.\\
(ii) Let $\varphi(x)=N_{F'/F}(N_C(x))$ with $N_C$ the quadratic norm of a composition algebra over $F'$,
$F'$ a quadratic field extension of $F$. $\varphi$ is a form of degree 4 permitting composition.
If $C$ has dimension greater than 1 then $\varphi$ satisfies SNP for all field extensions of odd degree (Lemma 3 (ii)).
If $C$ has dimension 1 then $\varphi$ satisfies SNP for all finite
field extensions  (Lemma 2 (i)).
For any form of degree 4 permitting composition, SNP holds for all odd degree Galois extensions.
\end{remark}

We conclude pointing out that already for cubic forms (which do not permit composition), it might not be enough any more to investigate if
$a\varphi_K\cong\varphi_K$ implies that $N_{K/F}(a)\varphi\cong\varphi$. It might also be interesting to know if and when
$N_{K/F}(a)^2\varphi\cong\varphi$ holds.

\smallskip
{\it Acknowledgments:} The second author would like to acknowledge the support of the
``Georg-Thieme-Ged\"{a}chtnisstiftung'' (Deutsche Forschungsgemeinschaft)
 during her stay at the University of Trento,
 and thanks the Department of Mathematics at Trento for its hospitality.
 We thank Barquero for pointing out [B1, B2] and [B-M] to us.

\end{document}